\documentclass[12pt]{article}
\usepackage{amsmath}
\usepackage{amssymb}
\numberwithin{equation}{section}

\title{Comment on the sums
$S(n)=\sum\limits_{k=-\infty}^\infty  \frac{1}{(4k+1)^n}$}

\author{ Z.~K.~Silagadze \\
Budker Institute of Nuclear Physics SB RAS and \\
Novosibirsk State University, 630 090, Novosibirsk, Russia }
\date{}

\begin{document}

\maketitle

\begin{abstract}
This is a comment on the papers N. D. Elkies, Amer. Math. Monthly  {\bf 110}  
(2003), 561--573 and  Cvijovi\'{c} and  J. Klinowski, J. Comput. Appl. Math.  
{\bf 142}  (2002), 435--439. We provide an explicit expression for
the kernel of the integral operator introduced in the first paper. This
explicit expression considerably simplifies the calculation of $S(n)$ and
enables a simple derivation of Cvijovi\'{c} and Klinowski's integral
representation for $\zeta(2n+1)$.
\end{abstract}

\section{Introduction}
In the remarkable paper \cite{1} Beukers, Kolk and Calabi showed that for all
$n\ge 2$,
\begin{equation}
S(n)=\sum\limits_{k=-\infty}^\infty \frac{1}{(4k+1)^n}=\left (\frac{\pi}{2}
\right )^n\,\delta_n,
\label{eq1}
\end{equation}
where $\delta_n=\mathrm{Vol}(\Delta_n)$ is the volume of the $n$-dimensional
convex polytope
\begin{equation}
\Delta_n=\left \{(u_1,\ldots,u_n): u_i > 0,\; u_i+u_{i+1} < 1
\right \}.
\label{eq2}
\end{equation}
Here $u_i$ are indexed cyclically modulo $n$, that is $u_{n+1}=u_1$.

In a subsequent paper \cite{2}, Elkies gave an elegant method for
calculating the volume $\delta_n$ (earlier calculations of this type can be
found in \cite{3}). If we introduce the characteristic function $K_1(u,v)$
of the isosceles right triangle $\{(u,v): u,v > 0,\,u+v <
1\}$ that is 1 inside  the triangle and 0 outside of it, then \cite{2}
\begin{eqnarray} &&
\delta_{n}=\int_0^1\ldots\int_0^1
\prod\limits_{i=1}^nK_1(u_i,u_{i+1})\,du_1\ldots du_n=\int_0^1du_1
\int_0^1du_2 \,K_1(u_1,u_2) \nonumber \\ &&
\ldots\int_0^1du_{n-1}\,K_1(u_{n-2},u_{n-1})\int_0^1du_n\,
K_1(u_{n-1},u_n)\,K_1(u_n,u_1).
\label{eq4}
\end{eqnarray}
Let us note that $\prod\limits_{i=1}^nK_1(u_i,u_{i+1})$ is the characteristic
function of the polytope $\Delta_n$ and just this property enables us to
extend the integration domain from $\Delta_n$ to the $n$-dimensional unit
hypercube and obtain (\ref{eq4}). At that the value of $\prod\limits_{i=1}^n
K_1(u_i,u_{i+1})$ on the boundary of $\Delta_n$ is, in fact, irrelevant and
in the above integral we can assume
\begin{equation}
K_1(u,v)=\theta(1-u-v),
\label{eq3}
\end{equation}
where $\theta(x)$ is the Heaviside step function with $\theta(0)=1/2$
(it doesn't matter, as far as the integral (\ref{eq4}) is concerned,  what
value is used for $\theta(0)$, the choice $\theta(0)=1/2$ will be convenient
for us in the following).

$K_1(u,v)$ can be interpreted \cite{2} as the kernel of the
linear operator $\hat T$ on the Hilbert space $L^2(0,1)$, defined as follows
\begin{equation}
(\hat T f)(u)=\int_0^1 K_1(u,v)f(v)\,dv=\int_0^{1-u} f(v)\,dv.
\label{eq5}
\end{equation}
Then (\ref{eq4}) shows that $\delta_{n}$ equals just the trace of the
operator $\hat T^n$:
\begin{equation}
\delta_{n}=\int_0^1 K_n(u_1,u_1)\,du_1.
\label{eq6}
\end{equation}
The  kernel $K_n(u,v)$ of $\hat T^n$ obeys the recurrence relation
\begin{equation}
K_n(u,v)=\int_0^1K_1(u,u_1)\,K_{n-1}(u_1,v)\,du_1.
\label{eq7}
\end{equation}

\section{An explicit expression for the kernel}
No solution of the recurrence relation (\ref{eq7}) was given in \cite{2}.
It is the aim of this short note to provide the solution. When $K_1(u,v)=
\theta(1-u-v)$, it is given by
\begin{eqnarray} &&
K_{2n}(u,v)=(-1)^n\,\frac{2^{2n-2}}{(2n-1)!}  \nonumber \\ &&
\times \left\{ \left [  E_{2n-1}\left (
\frac{u+v}{2}\right )+E_{2n-1}\left (\frac{u-v}{2}\right )\right ]
\theta(u-v)\right . \nonumber \\ &&
\left . +
\left [  E_{2n-1}\left (\frac{u+v}{2}\right )+
E_{2n-1}\left (\frac{v-u}{2}\right )\right ]\theta(v-u)\right\},
\label{eq8}
\end{eqnarray}
and
\begin{eqnarray} &&
K_{2n+1}(u,v)=(-1)^n\,\frac{2^{2n-1}}{(2n)!}  \nonumber \\ &&
\times \left\{ \left [  E_{2n}\left (
\frac{1-u+v}{2}\right )+E_{2n}\left (\frac{1-u-v}{2}\right )\right ]
\theta(1-u-v)\right . \nonumber \\ &&
 \left . +
\left [  E_{2n}\left (\frac{1-u+v}{2}\right )-
E_{2n}\left (\frac{u+v-1}{2}\right )\right ]\theta(u+v-1)\right\}.
\label{eq9}
\end{eqnarray}
In these formulas, $E_n(x)$ are the Euler polynomials (see e.g. \cite{4},
\cite{4a}) and $\theta(x)$ is the Heaviside step function with
$\theta(0)=1/2$.

After they have been guessed, it is quite straightforward  to prove by
induction that (\ref{eq8}) and (\ref{eq9}) obey the recurrence relation
(\ref{eq7}) (see the Appendix).

No systematic method was used to get (\ref{eq8}) and (\ref{eq9}). They have
indeed been guessed. We simply calculated a number of explicit expressions for
$K_n(u,v)$ using the recurrence relation (\ref{eq7}) and the expression
(\ref{eq3}), and then tried to locate regularities in these expressions.
The appearance of Euler polynomials suggests that $K_n(u,v)$ should have
a nice two-variable Fourier decomposition and a likely systematic approach
could be based on this fact \cite{2a}.

\section{Applications of the kernel}
Having (\ref{eq8}) and (\ref{eq9}) at our disposal, it is easy to calculate
the integral in (\ref{eq6}). Namely, because
\begin{equation}
K_{2n}(u,u)=(-1)^n\,\frac{2^{2n-2}}{(2n-1)!}\left [ E_{2n-1}(u)+ E_{2n-1}(0)
\right ],
\label{eq11}
\end{equation}
and
\begin{equation}
E_{2n-1}(u)=\frac{1}{2n}\,\frac{d}{du}\,E_{2n}(u),
\label{eq12}
\end{equation}
we get
\begin{equation}
\delta_{2n}=\int_0^1 K_{2n}(u,u)\,du=
(-1)^n\,\frac{2^{2n-2}}{(2n-1)!}E_{2n-1}(0)
\label{eq13}
\end{equation}
(note that $E_{2n}(0)=E_{2n}(1)=0$).
But $E_{2n-1}(0)$ can be expressed through the Bernoulli numbers
\begin{equation}
E_{2n-1}(0)=-\frac{2}{2n}\,(2^{2n}-1)B_{2n},
\label{eq14}
\end{equation}
and using the relation $S(2n)=(1-2^{-2n})\,\zeta(2n)$, we reproduce the
celebrated formula for the Riemann zeta function at even-integer arguments
\begin{equation}
\zeta(2n)=(-1)^{n+1}\,\frac{2^{2n-1}}{(2n)!}\,\pi^{2n}\,B_{2n}\,.
\label{eq15}
\end{equation}
In the case of $S(2n+1)$, we can use the identity
\begin{equation}
K_{2n+1}(u,u)+K_{2n+1}(1-u,1-u)=(-1)^n\;\frac{2^{2n}}{(2n)!}\;E_{2n}(1/2),
\label{eq16}
\end{equation} 
which follows from (\ref{eq9}), and get
\begin{eqnarray} &&
\delta_{2n+1} =\frac{1}{2}\int\limits_0^1
\left [ K_{2n+1}(u,u)+K_{2n+1}(1-u,1-u)\right ] du \nonumber \\ &&
=(-1)^n\,\frac{2^{2n-1}}{(2n)!}\;
E_{2n}(1/2).
\label{eq17}
\end{eqnarray}
Expressing  $E_{2n}(1/2)$ in terms of the Euler numbers
\[
E_{2n}(1/2)=2^{-2n}E_{2n},
\]
we get from (\ref{eq17})  Euler's other celebrated formula (note that
formulas (9) and (10) in \cite{2}, analogs of our (\ref{eq15}) and
(\ref{eq18}), contain typos)
\begin{equation}
S(2n+1)=(-1)^n\;\frac{1}{2(2n)!}\;\left (\frac{\pi}{2}\right )^{2n+1}E_{2n}.
\label{eq18}
\end{equation}
Interestingly, our expression for $K_n(u,v)$ allows us to re-derive
Cvijovi\'{c} and Klinowski's integral representation \cite{5} for
$\zeta(2n+1)$. We begin with the formula (for more details, see \cite{6})
\begin{equation}
\zeta(2n+1)=-\frac{2^{2n+1}}{2^{2n+1}-1}\,\frac{1}{2n}\,\idotsint\limits_
{\Box_{2n}}
\frac{\ln{(x_1\cdots x_{2n})}}{1-x^2_1\cdots x^2_{2n}}\,dx_1\cdots dx_{2n},
\label{eq19}
\end{equation}
where $\Box_{2n}$ is the $2n$-dimensional unit hypercube. If we now apply
the Beukers--Kolk--Calabi change of variables \cite{1}
\begin{equation}
x_1=\frac{\sin{u_1}}{\cos{u_2}},\;\;x_2=\frac{\sin{u_2}}{\cos{u_3}},\ldots,
\;x_{2n-1}=\frac{\sin{u_{2n-1}}}{\cos{u_{2n}}},\;\;x_{2n}=\frac{\sin{u_{2n}}}
{\cos{u_1}}
\label{eq20}
\end{equation}
to the integral (\ref{eq19}), after some simple manipulations we get
\begin{equation}
\zeta(2n+1)=-\frac{2^{2n+1}}{2^{2n+1}-1}\,\left(\frac{\pi}{2}\right)^{2n}
\idotsint\limits_{\Delta_{2n}}\ln{\left[\tan{\left(u_1\frac{\pi}{2}\right)}
\right]}\,du_1\cdots du_{2n}.
\label{eq21}
\end{equation}
Using the kernel $K_{2n}(u,v)$, we can reduce the evaluation of (\ref{eq21})
to the evaluation of the following one-dimensional integral:
\begin{equation}
\zeta(2n+1)=-\frac{2\,\pi^{2n}}{2^{2n+1}-1}\int_0^1
\ln{\left[\tan{\left(\frac{\pi}{2}u\right)}\right]}\,K_{2n}(u,u)\,du.
\label{eq22}
\end{equation}
But
\[
\ln{\left[\tan{\left(\frac{\pi}{2}(1-u)\right)}\right]}=
\ln{\left[\cot{\left(\frac{\pi}{2}u\right)}\right]}=
-\ln{\left[\tan{\left(\frac{\pi}{2}u\right)}\right]},
\]
which enables us to rewrite (\ref{eq22}) as
\begin{eqnarray} &&
\zeta(2n+1) =-\frac{\pi^{2n}}{2^{2n+1}-1} \nonumber \\ &&
\times \int_0^1
\ln{\left[\tan{\left(\frac{\pi}{2}u\right)}\right]}\,\left [
K_{2n}(u,u)- K_{2n}(1-u,1-u)\right ]\,du.
\label{eq23}
\end{eqnarray}
However, from (\ref{eq11}) and (\ref{eq12}) we have (recall that
$E_{2n-1}(1-u)=-E_{2n-1}(u))$
\[
K_{2n}(u,u)- K_{2n}(1-u,1-u)=(-1)^n\,\frac{2^{2n-1}}{(2n)!}\,
\frac{d}{du}E_{2n}(u),
\]
and the straightforward integration by parts in (\ref{eq23}) yields finally
the result
\begin{equation}
\zeta(2n+1)=\frac{(-1)^n\,\pi^{2n+1}}{4\,[1-2^{-(2n+1)}]\,(2n)!}\,
\int_0^1\frac{E_{2n}(u)}{\sin{(\pi\,u)}}\,du.
\label{eq24}
\end{equation}
This is exactly the integral representation for $\zeta(2n+1)$ found in
\cite{5}.

\section{Concluding remarks}
Finally, let us comment on the origin of the highly non-trivial
Beukers--Kolk--Calabi change of variables. Using the hyperbolic version of it
\begin{equation}
x_1=\frac{\sinh{u_1}}{\cosh{u_2}},\;\;x_2=\frac{\sinh{u_2}}{\cosh{u_3}},
\ldots, x_{n-1}=\frac{\sinh{u_{n-1}}}{\cosh{u_n}},\;\;x_n=\frac{\sinh{u_n}}
{\cosh{u_1}},
\label{eq25}
\end{equation}
we can get \cite{6}
\begin{equation}
\zeta(n)=\frac{2^n}{2^n-1}\idotsint\limits_{U_n}du_1\cdots du_n=
\frac{2^n}{2^n-1}\;\mathrm{Vol}_{n}(U_n),
\label{eq26}
\end{equation}
where $U_n$ has a complicated amoeba-like shape with narrowing tentacles
going to infinity. It was shown by Passare \cite{7} that $\zeta(2)$ is
really related to the amoeba, a fascinating object in complex geometry
(\cite{8}, \cite{9}) introduced in the book \cite{10}. We can conjecture that
$U_n$ is the $2^n$-th part of the amoeba  of a certain Laurent polynomial and
the Newton polytope of this polynomial is the union, from which the central
point $(0,\ldots,0)$ is thrown away, of the (suitably rescaled) polytope
$\Delta_n$ and its $2^n-1$ mirror images under reflections $u_i\to -u_i$.
For $n=2$, the conjecture is valid and the corresponding amoeba uncovers the
origin of the two-dimensional version of the Beukers--Kolk--Calabi change of
variables \cite{6}. We think it is worthwhile to further investigate  this
possible connection between $\zeta(n)$ and amoebas.

\appendix
\section{Appendix: The proof of (\ref{eq8}) and (\ref{eq9}) by induction}
If $n=0$, (\ref{eq9}) gives, because of $E_0(x)=1$,
\[
K_1(u,v)=\theta(1-u-v),
\]
as desired. Let us calculate $K_2(u,v)$. From (\ref{eq7}) we have
$$K_2(u,v)=\int\limits_0^1K_1(u,u_1)\,K_1(u_1,v)\,du_1= $$ $$=
\int\limits_0^{1-u}\theta(1-v-u_1)\,du_1=\left\{\begin{array}{c}
1-v,\;\;\mathrm{if}\;\;v\ge u, \\ 1-u,\;\;\mathrm{if}\;\;v\le u.
\end{array}\right . $$
We can rewrite this as follows
\[
K_2(u,v)=(1-v)\theta(v-u)+(1-u)\theta(u-v),
\]
and this is exactly what (\ref{eq3}) gives for $n=1$ because $E_1(x)=x-1/2$.

Now suppose (\ref{eq8}) is true for a particular $n$ and let us calculate
\[
K_{2n+1}(u,v)=\int_0^1K_1(u,u_1)\,K_{2n}(u_1,v)\,du_1.
\]
We will use the  following property of the Euler polynomials
\begin{equation}
\frac{d}{dx}E_n(x)=nE_{n-1}(x),
\label{eqa5}
\end{equation}
which facilitates the calculation of integrals.

We have
\begin{eqnarray} &&
\int_0^{1-u}E_{2n-1}\left(\frac{u_1+v}{2}\right )\,\theta(u_1-v)\,
du_1 \nonumber \\ &&
=\theta(1-u-v)\int_v^{1-u}E_{2n-1}\left(\frac{u_1+v}{2}\right )
\,du_1.
\end{eqnarray}
Let us make the change of variables
\[
\frac{u_1+v}{2}=x.
\]
Then the integral equals to
\[
2\int_v^{(1-u+v)/2}\frac{1}{2n}\left (\frac{d}{dx}E_{2n}(x)\right )
dx=\frac{1}{n}\left (E_{2n}\left(\frac{1-u+v}{2}\right)-E_{2n}(v)\right).
\]
Therefore,
\begin{eqnarray} &&
\int_0^{1-u}E_{2n-1}\left(\frac{u_1+v}{2}\right )\,\theta(u_1-v)\,
du_1 \nonumber \\ &&
=\frac{\theta(1-u-v)}{n}\left [E_{2n}\left(\frac{1-u+v}{2}\right)-
E_{2n}(v)\right].
\label{eqa6}
\end{eqnarray}
Analogously, using the substitution $\frac{u_1-v}{2}=x$ and $E_{2n}(0)=0$,
we get
\begin{equation}
\int_0^{1-u}E_{2n-1}\left(\frac{u_1-v}{2}\right )\,\theta(u_1-v)\,
du_1=\frac{\theta(1\!-\!u\!-\!v)}{n}E_{2n}\left(\frac{1-u-v}{2}\right).
\label{eqa7}
\end{equation}
Further we have
\begin{eqnarray} &&
\int_0^{1-u}E_{2n-1}\left(\frac{u_1+v}{2}\right ) \theta(v-u_1)
du_1 \nonumber \\ &&
=\theta(1\!-\!u\!-\!v)\int_0^v \!E_{2n-1}\left(\frac{u_1\!+\!v}{2}
\right )du_1 \nonumber \\ &&
+\theta(u\!+\!v\!-\!1)\int_0^{1-u} \!E_{2n-1}
\left(\frac{u_1\!+\!v}{2}\right )du_1 \nonumber \\ &&
=\frac{\theta(1-u-v)}{n}\left[E_{2n}(v)-E_{2n}\left(\frac{v}{2}\right)
\right]
\nonumber \\ &&
+\frac{\theta(u+v-1)}{n}\left[E_{2n}\left(\frac{1-u+v}{2}\right)-
E_{2n}\left(\frac{v}{2}\right)\right],
\label{eqa8}
\end{eqnarray}
and analogously
\begin{eqnarray} &&
\int_0^{1-u}E_{2n-1}\left(\frac{v-u_1}{2}\right )\,\theta(v-u_1)\,
du_1=
\frac{\theta(1-u-v)}{n}\,E_{2n}\left(\frac{v}{2}\right) \nonumber \\ &&
+
\frac{\theta(u+v-1)}{n}\left[E_{2n}\left(\frac{v}{2}\right)-
E_{2n}\left(\frac{u+v-1}{2}\right)\right].
\label{eqa9}
\end{eqnarray}
The equations (\ref{eq8}), (\ref{eqa6}), (\ref{eqa7}), (\ref{eqa8}) and
(\ref{eqa9}) indicate that $K_{2n+1}(u,v)$ is given by (\ref{eq9}), as
desired.

Now let us show that the validity of (\ref{eq9}) implies that
\begin{equation}
K_{2n+2}(u,v)=\int_0^1K_1(u,u_1)\,K_{2n+1}(u_1,v)\,du_1
\label{eqa10}
\end{equation}
is given by the formula (\ref{eq8}).

We have
\begin{eqnarray} &&
\int_0^{1-u}E_{2n}\left (\frac{1-u_1+v}{2}\right )\theta(1-u_1-v)\,
du_1 \nonumber \\ &&
=\theta(v-u)\int_0^{1-v}E_{2n}\left (\frac{1-u_1+v}{2}
\right )du_1 \nonumber \\ &&
+\theta(u-v)\int_0^{1-u}E_{2n}\left (\frac{1-u_1+v}{2}
\right )du_1 = \nonumber \\ &&
-\frac{2}{2n+1}\,\theta(v-u)\left [E_{2n+1}(v)-E_{2n+1}
\left (\frac{1+v}{2}\right)\right ] \nonumber \\ &&
-\frac{2}{2n+1}\,\theta(u-v)\left [E_{2n+1}\left (\frac{u+v}{2}\right )-
E_{2n+1}\left (\frac{1+v}{2}\right)\right ],
\label{eqa11}
\end{eqnarray}
and
\begin{eqnarray} &&
\int_0^{1-u}E_{2n}\left (\frac{1-u_1-v}{2}\right )\theta(1-u_1-v)\,du_1
\nonumber \\ &&
=-\frac{2}{2n+1}\,\theta(v-u)\left [E_{2n+1}(0)-
E_{2n+1}\left (\frac{1-v}{2}\right)\right ]
\nonumber \\ &&
-\frac{2}{2n+1}\,\theta(u-v)\left [E_{2n+1}\left (\frac{u-v}{2}\right )-
E_{2n+1}\left (\frac{1-v}{2}\right)\right ],
\label{eqa12}
\end{eqnarray}
as well as
\begin{eqnarray} &&
\int_0^{1-u}E_{2n}\left (\frac{1-u_1+v}{2}\right )\theta(u_1+v-1)\, du_1 
\nonumber \\ &&
=-\frac{2}{2n+1}\,\theta(v-u)\left [E_{2n+1}\left (\frac{u+v}{2}
\right )-E_{2n+1}(v)\right ],
\label{eqa13}
\end{eqnarray}
and finally
\begin{eqnarray} &&
\int_0^{1-u}E_{2n}\left (\frac{u_1+v-1}{2}\right )\theta(u_1+v-1)\, du_1 
\nonumber \\ &&
=\frac{2}{2n+1}\,\theta(v-u)\left [E_{2n+1}\left (\frac{v-u}{2}
\right )-E_{2n+1}(0)\right ].
\label{eqa14}
\end{eqnarray}
In light of (\ref{eq9}), (\ref{eqa11}), (\ref{eqa12}), (\ref{eqa13}) and
(\ref{eqa14}), the integral (\ref{eqa10}) takes the form of (\ref{eq8}),
with $n\to n+1$, because
\begin{equation}
E_{2n+1}\left (\frac{1-v}{2} \right )+E_{2n+1}\left (\frac{1+v}{2}\right )=0,
\label{eqa16}
\end{equation}
which follows from $E_n(1-x)=(-1)^nE_n(x)$.

\section*{Acknowledgments}
The work is supported by the Ministry of Education and
Science of the Russian Federation and in part by Russian Federation President
Grant for the support of scientific schools NSh-5320.2012.2.


\begin{thebibliography}{99}
\bibitem{1}
F. Beukers, J.~A.~C. Kolk,  and E. Calabi,
Sums of generalized harmonic series and volumes,
\emph{Nieuw Arch. Wisk. (4)}  \textbf{11}:3  (1993), 217--224.

\bibitem{4a}
G. Bretti and P. E.  Ricci,
Euler polynomials and the related quadrature rule,
\emph{Georgian Math. J.}  \textbf{8}:3  (2001), 447--453.

\bibitem{5}
D. Cvijovi\'{c} and  J. Klinowski,
Integral representations of the Riemann zeta function for odd-integer
arguments,
\emph{J. Comput. Appl. Math. } \textbf{142}:2  (2002), 435--439.

\bibitem{2}
N. D. Elkies,
On the sums $\sum\sp \infty\sb {k=-\infty}(4k+1)\sp {-n}$,
\emph{Amer. Math. Monthly } \textbf{110}:7  (2003), 561--573;
Corr. ibid. {\bf 111} (2004), 456.

\bibitem{2a}
N. D. Elkies, Personal communication, 2010.

\bibitem{10}
I. M. Gel'fand,   M. M. Kapranov,  and A. V. Zelevinsky,
\emph{Discriminants, Resultants, and Multidimensional Determinants},
Mathematics: Theory \& Applications. Birkh\"auser Boston, Inc., Boston, MA,
1994.

\bibitem{4}
B. K. Karande and N. K. Thakare,
On the unification of Bernoulli and Euler polynomials,
\emph{Indian J. Pure Appl. Math.}  \textbf{6}:1  (1975), 98--107.

\bibitem{3}
J. Kubilius,
Estimation of the second central moment for strongly additive arithmetic
functions,
\emph{Litovsk. Mat. Sb.}  \textbf{23}:1  (1983), 122--133 (in Russian).

\bibitem{7}
M. Passare,
How to compute $\sum 1/n\sp 2$ by solving triangles,
\emph{Amer. Math. Monthly}  \textbf{115}:8  (2008), 745--752.

\bibitem{9}
M. Passare  and A. Tsikh,
Amoebas: their spines and their contours,
in: \emph{Idempotent mathematics and mathematical physics},  pp. 275--288,
Contemp. Math., 377, Amer. Math. Soc., Providence, RI, 2005.

\bibitem{6}
Z. K. Silagadze,
\emph{Sums of generalized harmonic series for kids from five to fifteen},
preprint, 2010; arXiv:1003.3602 [math.CA] (2010).

\bibitem{8}
O. Viro,
What 1s an amoeba?
\emph{Notices Amer.\ Math.\ Soc.}  {\bf 49}:8 (2002), 916--917. 

\end{thebibliography}
\end{document}